\magnification=\magstep1


\def\item{\vskip1.3pt\hang\textindent}


\tolerance=300 \pretolerance=200 \hfuzz=1pt \vfuzz=1pt

\hoffset 0cm            
\hsize=5.8 true in \vsize=9.5 true in

\def\rightheadline{\hfil\smc\lastname\hfil\tenbf\folio}
\def\leftheadline{\tenbf\folio\hfil\smc\lastname\hfil}
\headline={\ifodd\pageno\rightheadline\else\leftheadline\fi}
\newdimen\dimenone
\def\checkleftspace#1#2#3#4#5{
 \dimenone=\pagetotal
 \advance\dimenone by -\pageshrink   
 \ifdim\dimenone>\pagegoal          
   \else\dimenone=\pagetotal
        \advance\dimenone by \pagestretch
        \ifdim\dimenone<\pagegoal
          \dimenone=\pagetotal
          \advance\dimenone by#1         
          \setbox0=\vbox{#2\parskip=0pt                
                       \hyphenpenalty=10000
                       \rightskip=0pt plus 5em
                       \noindent#3 \vskip#4}    
        \advance\dimenone by\ht0
        \advance\dimenone by 3\baselineskip
        \ifdim\dimenone>\pagegoal\vfill\eject\fi
          \else\eject\fi\fi}

\parindent=35pt
\mathsurround=1pt
\parskip=1pt plus .25pt minus .25pt
\normallineskiplimit=.99pt

\mathchardef\emptyset="001F 

\def\summe#1{\displaystyle\mathop{\textstyle\sum^{#1}}\limits}
\def\summeint#1#2{\displaystyle\mathop{\textstyle\sum^{#1} \int \d{#2}}\limits}
\def\summedoint#1#2{\displaystyle\mathop{\textstyle\sum^{#1} \int\int\d{#2}}\limits}

\def\Int{\mathop{\rm int}\nolimits}
%



\def\1{{\bf1}}\def\0{{\bf0}}

\def\({\bigl(}  \def\){\bigr)}
\def\<{\mathopen{\langle}}\def\>{\mathclose{\rangle}}

\def\Z{{\mathchoice{{\hbox{$\rm Z\hskip 0.26em\llap{\rm Z}$}}}%
{{\hbox{$\rm Z\hskip 0.26em\llap{\rm Z}$}}}%
{{\hbox{$\scriptstyle\rm Z\hskip 0.31em\llap{$\scriptstyle\rm Z$}$}}}{{%
\hbox{$\scriptscriptstyle\rm
Z$\hskip0.18em\llap{$\scriptscriptstyle\rm Z$}}}}}}

\def\R{{\mathchoice{\hbox{$\rm I\hskip-0.14em R$}}%
{\hbox{$\rm I\hskip-0.14em R$}}%
{\hbox{$\scriptstyle\rm I\hskip-0.14em R$}}%
{\hbox{$\scriptscriptstyle\rm I\hskip-0.10em R$}}}}

\def\.{{\cdot}}
\def\|{\Vert}
\def\ssk{\smallskip}
\def\msk{\medskip}
\def\bsk{\bigskip}
\def\giantskip{\vskip2\bigskipamount}

\def\giantbreak{\par \ifdim\lastskip<2\bigskipamount \removelastskip
         \penalty-400 \giantskip\fi}

\def\nin{\noindent}
\def\cen{\centerline}
\def\pagebreak{\vskip 0pt plus 0.0001fil\break}
\def\linebreak{\break}

\def\epsilon{\varepsilon}

\font\ninerm=cmr9 \font\eightrm=cmr8 \font\sixrm=cmr6
 \font\eightbf=cmbx8 \font\sixbf=cmbx6
 \font\eighti=cmmi8 \font\sixi=cmmi6
\font\ninesy=cmsy9 \font\eightsy=cmsy8 \font\sixsy=cmsy6
 \font\eightit=cmti8 
 \font\eightsl=cmsl8 
\font\eighttt=cmtt8
\font\bfone=cmbx10 scaled\magstep1 
\font\smc=cmcsc10 
 
scaled\magstep1 \font\small=cmcsc8

\def\no #1. {\bigbreak\vskip-\parskip\noindent\bf #1. \quad\rm}

\def\Proposition #1. {\checkleftspace{0pt}{\bf}{Theorem}{0pt}{}
\bigbreak\vskip-\parskip\noindent{\bf Proposition #1.} \quad\it}

\def\Theorem #1. {\checkleftspace{0pt}{\bf}{Theorem}{0pt}{}
\bigbreak\vskip-\parskip\noindent{\bf  Theorem #1.} \quad\it}
\def\Corollary #1. {\checkleftspace{0pt}{\bf}{Theorem}{0pt}{}
\bigbreak\vskip-\parskip\nin{\bf Corollary #1.} \quad\it}
\def\Lemma #1. {\checkleftspace{0pt}{\bf}{Theorem}{0pt}{}
\bigbreak\vskip-\parskip\noindent{\bf  Lemma #1.}\quad\it}

\def\Definition #1. {\checkleftspace{0pt}{\bf}{Theorem}{0pt}{}
\rm\bigbreak\vskip-\parskip\noindent{\bf Definition #1.} \quad}

\def\Remark #1. {\checkleftspace{0pt}{\bf}{Theorem}{0pt}{}
\rm\bigbreak\vskip-\parskip\noindent{\bf Remark #1.}\quad}

\def\Exercise #1. {\checkleftspace{0pt}{\bf}{Theorem}{0pt}{}
\rm\bigbreak\vskip-\parskip\noindent{\bf Exercise #1.} \quad}

\def\Example #1. {\checkleftspace{0pt}{\bf}{Theorem}{0pt}{}
\rm\bigbreak\vskip-\parskip\noindent{\bf Example #1.}\quad}
\def\Examples #1. {\checkleftspace{0pt}{\bf}{Theorem}{0pt}
\rm\bigbreak\vskip-\parskip\noindent{\bf Examples #1.}\quad}

\newcount\problemnumb \problemnumb=0
\def\Problem{\global\advance\problemnumb by 1\bigbreak\vskip-\parskip\noindent
{\bf Problem \the\problemnumb.}\quad\rm }

\def\Proof#1.{\rm\par\ifdim\lastskip<\bigskipamount\removelastskip\fi\smallskip
            \noindent {\bf Proof.}\quad}

\nopagenumbers

\def\author{}
\def\lastname{}
\def\thanks#1{\footnote*{\eightrm#1}}
\def\title{}

\def\lastname{}
\def\h{{\textstyle{1\over2}}}

\def\he{{1\over2}}

\def\ep{\epsilon}

\def\text{\textstyle}
\def\disp{\displaystyle}
\def\d{{\,\rm d}}

\def\and{{\rm and }}

\expandafter\edef\csname amssym.def\endcsname{%
       \catcode`\noexpand\@=\the\catcode`\@\space}
\catcode`\@=11
\def\undefine#1{\let#1\undefined}
\def\newsymbol#1#2#3#4#5{\let\next@\relax
 \ifnum#2=\@ne\let\next@\msafam@\else
 \ifnum#2=\tw@\let\next@\msbfam@\fi\fi
 \mathchardef#1="#3\next@#4#5}
\def\mathhexbox@#1#2#3{\relax
 \ifmmode\mathpalette{}{\m@th\mathchar"#1#2#3}%
 \else\leavevmode\hbox{$\m@th\mathchar"#1#2#3$}\fi}
\def\hexnumber@#1{\ifcase#1 0\or 1\or 2\or 3\or 4\or 5\or 6\or 7\or 8\or
 9\or A\or B\or C\or D\or E\or F\fi}


\font\tenmsb=msbm10 \font\sevenmsb=msbm7 \font\fivemsb=msbm5
\newfam\msbfam
\textfont\msbfam=\tenmsb\scriptfont\msbfam=\sevenmsb
\scriptscriptfont\msbfam=\fivemsb \edef\msbfam@{\hexnumber@\msbfam}
\def\Bbb#1{{\fam\msbfam\relax#1}}

\font\teneufm=eufm10 \font\seveneufm=eufm7 \font\fiveeufm=eufm5
\newfam\eufmfam
\textfont\eufmfam=\teneufm \scriptfont\eufmfam=\seveneufm
\scriptscriptfont\eufmfam=\fiveeufm

\catcode`@=11 

\expandafter\edef\csname amssym.def\endcsname{%
       \catcode`\noexpand\@=\the\catcode`\@\space}
\font\eightmsb=msbm8 \font\sixmsb=msbm6 \font\fivemsb=msbm5
\font\eighteufm=eufm8 \font\sixeufm=eufm6 \font\fiveeufm=eufm5
\newskip\ttglue
\def\eightpoint{\def\rm{\fam0\eightrm}%
  \textfont0=\eightrm \scriptfont0=\sixrm \scriptscriptfont0=\fiverm
  \textfont1=\eighti \scriptfont1=\sixi \scriptscriptfont1=\fivei
  \textfont2=\eightsy \scriptfont2=\sixsy \scriptscriptfont2=\fivesy
  \textfont3=\tenex \scriptfont3=\tenex \scriptscriptfont3=\tenex
\textfont\eufmfam=\eighteufm \scriptfont\eufmfam=\sixeufm
\scriptscriptfont\eufmfam=\fiveeufm \textfont\msbfam=\eightmsb
\scriptfont\msbfam=\sixmsb \scriptscriptfont\msbfam=\fivemsb
  \def\it{\fam\itfam\eightit}%
  \textfont\itfam=\eightit
  \def\sl{\fam\slfam\eightsl}%
  \textfont\slfam=\eightsl
  \def\bf{\fam\bffam\eightbf}%
  \textfont\bffam=\eightbf \scriptfont\bffam=\sixbf
   \scriptscriptfont\bffam=\fivebf
  \def\tt{\fam\ttfam\eighttt}%
  \textfont\ttfam=\eighttt
  \tt \ttglue=.5em plus.25em minus.15em
  \normalbaselineskip=9pt
  \def\MF{{\manual opqr}\-{\manual stuq}}%
  \let\big=\eightbig
  \setbox\strutbox=\hbox{\vrule height7pt depth2pt width\z@}%
  \normalbaselines\rm}
\def\eightbig#1{{\hbox{$\textfont0=\ninerm\textfont2=\ninesy
  \left#1\vbox to6.5pt{}\right.\n@space$}}}


\csname amssym.def\endcsname


\def\la{\lambda}
\def\al{\alpha}
\def\be{\beta}

\def\({\left(}
\def\){\right)}
\def\for{\qquad \hbox{for}\ }
\def\eq{\eqalign}

\def\O#1{O\(#1\)}
\def\abs#1{\left| #1 \right|}

\def\norm#1{\left\Vert #1 \right\Vert}

\def\klein{\eightpoint \def\smc{\small} \baselineskip=9pt}

\def\fn#1#2{{\parindent=0.7true cm
\footnote{$^{(#1)}$}{{\klein  #2}}}}

\font\boldmas=msbm10                  
\def\Bbb#1{\hbox{\boldmas #1}}        
\def\Z{{\Bbb Z}}                        

\def\R{{\Bbb R}}


\font\eightrm=cmr8 \long\def\fussnote#1#2{{\baselineskip=9pt
\setbox\strutbox=\hbox{\vrule height 7pt depth 2pt width 0pt}%
\eightrm \footnote{#1}{#2}}}
\font\boldmasi=msbm10 scaled 700      
\def\Bbbi#1{\hbox{\boldmasi #1}}      
\font\boldmas=msbm10                  
\def\Bbb#1{\hbox{\boldmas #1}}        
\def\Zi{{\Bbbi Z}}                      
\def\Pi{{\Bbbi P}}                      
\def\Ri{{\Bbbi R}}



\def\dint #1 {
\quad  \setbox0=\hbox{$\disp\int\!\!\!\int$}
  \setbox1=\hbox{$\!\!\!_{#1}$}
  \vtop{\hsize=\wd1\centerline{\copy0}\copy1} \quad}

\def\drint #1 {
\qquad  \setbox0=\hbox{$\disp\int\!\!\!\int\!\!\!\int$}
  \setbox1=\hbox{$\!\!\!_{#1}$}
  \vtop{\hsize=\wd1\centerline{\copy0}\copy1}\qquad}

\def\frac#1#2{{#1\over #2}}

\def\date{\the\day.~\the\month.~\the\year}

\def\klein{\eightpoint \def\smc{\small} }

\def\at#1#2#3{{\left.\phantom{\Int} #1 \right|}_{#2}^{#3}}

\def\frac#1#2{{#1\over#2}}
\def\Int{\int\limits}

\def\vol{{\rm vol}}

\hsize=16true cm     \vsize=23true cm \parindent=0pt

\def\cS{{\cal S}} \def\cB{{\cal B}} \def\cC{{\cal C}}
\def\J{{\cal J}} \def\cL{{\cal L}}
\def\cR{{\cal R}} \def\cF{{\cal F}} \def\l{\ell} \def\k{q}
\def\cU{{\cal U}} \def\cV{{\cal V}} \def\cI{{\cal I}}

\cen{{\bfone The lattice discrepancy of certain three-dimensional
bodies}} \bsk 

\cen{{\bf Ekkehard Kr\"{a}tzel and Werner Georg Nowak }\fn{*}{The
authors gratefully acknowledge support from the Austrian Science
Fund (FWF) under project Nr.~P20847-N18.} {\bf(Vienna)} }

\vbox{\vskip 1.2true cm}

{\klein{\bf Abstract. } Based on a very precise approximation to
the lattice discrepancy of a Lam\'e disc, an asymptotic formula is
established for the number of lattice points in the
three-dimensional body $$ |u_1|^{mk}+ \(|u_2|^k+|u_3|^k\)^{m} \le
x^{mk}\,,$$ for large real $x$ and fixed reals $m, k$. Particular
attention is paid to the boundary points of Gaussian curvature
zero. \bsk

bodies, exponential sums}

\vbox{\vskip 1true cm}

{\bf 1.~Introduction. } For a compact body $\cB$ in $\R^3$ the
{\it lattice point discrepancy } of a copy $x\cB$, obtained by
linear dilation by a large real parameter $x$, is of interest:
$$ P_\cB(x) := \#\(x\cB\cap\Z^3\)-\vol(\cB)x^3\,. \eqno(1.1)  $$
The theory of its estimation, resp., asymptotic evaluation has
been described in detail in Kr\"{a}tzel's books [13], [15], even in
the frame of a more general $s$-dimensional setting, and also in
Huxley's monograph [7] which concentrates on the planar case. See
also the survey article [11] by the authors with Ivi\'c and
K\"{u}hleitner. \ssk The case of a convex body $\cC$ of smooth
boundary $\partial\cC$ with bounded nonzero Gaussian curvature
throughout is quite well understood. Hlawka's [6] classic bound
$P_\cC(x)\ll x^{3/2}$ has been refined up to $O(x^{63/43+\ep})$ by
W.~M\"{u}ller [23], using exponential sums in a most ingenious way.
Furthermore, we also know that
$$P_\cC(x)=\Omega\(x(\log x)^{1/3}\)\,,\qquad\and\qquad \Int_0^X
(P_\cC(x))^2\d x \ll X^3(\log X)^2\,. $$ See Nowak [24], and
Iosevich, Sawyer \& Seeger [9]. For the special case of the sphere
in $\R^3$, cf.~Heath-Brown [5], for rational ellipsoids Chamizo,
Crist\'obal \& Ubis [1], for ellipsoids of rotation the present
authors [19].\ssk As soon as the boundary $\partial\cB$ has points
of curvature zero, the problem becomes much harder and,
consequently, our knowledge is rather fragmentary. For partial
results, see Haberland [4], Kr\"{a}tzel [14], [16], [17], [18] and
Peter [26]. A fairly general theorem on {\it bodies of rotation }
has been established by Nowak [25]. The method used there can be
called the {\it"cut-into-slices approach"}: Since $x\cB$
intersects every plane orthogonal to its axis of rotation in a
circular disc, a very accurate approximation to the latter's
discrepancy is employed, followed by a ("careful") summation with
respect to the third coordinate. This idea has also been pursued
by the authors for the special case of a rotating Lam\'e's curve
[20].\msk As a rule, with this sort of problems one cannot expect
results which are both very sharp {\it and } general. Since the
methods used are highly technical, in order to get precise
estimates one can only deal with a particular class of bodies at a
time. For the expert it will be transparent that the argument
developed may be applied in other similar situations as well. \ssk
On the basis of this understanding, in the present paper we will
concentrate on bodies
$$ \cB_{m,k}:=\{(u_1,u_2,u_3)\in\R^3:\ |u_1|^{mk}+ \(|u_2|^k+|u_3|^k\)^{m} \le 1\
\}\,, \eqno(1.3)$$ where $k>2$ and $m>1$ are fixed real numbers.

\msk

Let $A_{m,k}(x)$ the number of integer points in $x\cB_{m,k}$,
then for integers $k$ and $m$, $A_{m,k}(x)$ describes the average
number of representations of integers by the form $|u_1|^{mk}+
\(|u_2|^k+|u_3|^k\)^{m}$: $$ A_{m,k}(T^{1/(mk)})=\sum_{n\le
T}r_{m,k}(n)\,,$$ where $$ r_{m,k}(n):=\#\{ (n_1,n_2,n_3)\in\Z^3:\
|n_1|^{mk}+ \(|n_2|^k+|n_3|^k\)^{m}=n \}\,. $$ The cases $m=1$
(super sphere) and $k=2$ (rotating Lam\'e disc) have been dealt
with by Kr\"{a}tzel [13], [14], resp., by the authors [20]. \ssk An
important issue of the present paper is that the cut-into-slices
method will be applied for the first time to a body which lacks
symmetry of rotation. Thus, the planar discs arising are no longer
circular, but they are bounded by Lam\'e curves. The latter also
contain points of curvature zero, which renders the situation much
more complicated. We will establish sort of a "truncated Hardy's
identity for Lam\'e discs" which should be of some interest of its
own - see Theorem 2 below. \ssk The body $x\cB_{m,k}$ has been
investigated earlier by Kr\"{a}tzel [18], by a purely
"three-dimensional" method: This lead to the estimation of double
trigonometric sums and to a somewhat less precise final result.
\ssk In fact, the Gaussian curvature of $\partial\cB_{m,k}$
vanishes on each curve of intersection with one of the coordinate
planes. Where any two of these curves meet, there arise {\it flat
points}, i.e., in $(\pm1,0,0)$, $(0,\pm1,0)$, and $(0,0,\pm1)$. We
will be able to evaluate precisely the contribution of these flat
points to the lattice discrepancy. The rest of the curves of
Gaussian curvature zero will be taken into account by appropriate
$O$-terms. (Cf.~the remarks after the statement of the main
theorem.)

\bsk\bsk

{\bf Theorem 1. } { \it For fixed reals $k>2$, $m>1$, with
$mk\ge{7\over3}$, and large real $x$, the number $A_{m,k}(x)$ of
lattice points in the body $x\cB_{m,k}$ satisfies the asymptotic
formula $$ A_{m,k}(x) = \vol(\cB_{m,k})x^3 + H_{mk,k,1}(x)+
H_{mk,k,2}(x) + \cR_{m,k}(x)\,, $$ where $$ \eq{ H_{mk,k,1}(x) &=
\cF_{mk,k,1}(x)\,x^{2-2/mk} + O(x)\,,\cr H_{mk,k,2}(x) &=
\cF_{mk,k,2}(x)\,x^{2-1/mk-1/k} + O(x)\,,\cr} $$ with continuous
periodic functions $\cF_{mk,k,i}(x)$, $i=1,2$. The functions $H$
and $\cF$ will be defined and analyzed in section 2. The remainder
term is given by
$$ \eq{\cR_{m,k}(x) = \O{x^{37/25}} &+
\O{x^{{339\over208}-{131\over104mk}}(\log
x)^{18627mk-20614\over8320mk}}\cr &+
\O{x^{{339\over208}-{235\over208k}}(\log
 x)^{{18627\over8320}(1-1/k)}}\,.\cr}  $$ Thus,} $$ \cR_{m,k}(x)
 \ll\ \cases{x^{37/25} & \hskip-3.5true cm if $k<{5875\over779}=7.54\dots$,
and $mk<{6550\over779}=8.408\dots$,\cr
x^{{339\over208}-{131\over104mk}}(\log x)^{18627mk-20614\over8320mk}
& if $mk\ge{6550\over779}$, and $m\ge{262\over235}=1.11\dots$, \cr
x^{{339\over208}-{235\over208k}}(\log
 x)^{{18627\over8320}(1-1/k)} & if $k\ge{5875\over779}$, and $m<{262\over235}$.\cr} $$

\bsk\bsk

{\bf Remarks. } As hinted at above, the term
$\cF_{mk,k,1}(x)\,x^{2-2/mk}$ comes from the points $(\pm1,0,0)$
where the two Lam\'e curves $$
\cL_1:\quad\cases{|u_1|^{mk}+|u_2|^{mk}=1,& \cr u_3=0\cr}
\qquad\and\qquad\cL_2:\quad\cases{|u_1|^{mk}+|u_3|^{mk}=1,&\cr
u_2=0\cr} $$ intersect. The curvature of either of these curves
vanishes at $(\pm1,0,0)$, of order $mk-2$ (related to the
arclength). Similarly, the points $(0,\pm1,0)$, $(0,0,\pm1)$ which
lie at the intersection of $\cL_1$, resp., $\cL_2$ with $$
\cL_3:\quad\cases{|u_2|^{k}+|u_3|^{k}=1,& \cr u_1=0,\cr}
$$ contribute the term $\cF_{mk,k,2}(x)\,x^{2-1/mk-1/k}$.
The curvature of $\cL_3$ has zeros of order $k-2$ at $(0,\pm1,0)$,
$(0,0,\pm1)$. As the proof will show, the contribution of the rest
of $\cL_1$, $\cL_2$ to the lattice point discrepancy can be bounded
by $ \O{x^{{339\over208}-{131\over104mk}}(\log
x)^{18627mk-20614\over8320mk}} $ while $\cL_3$ without the flat
points contributes at most $
\O{x^{{339\over208}-{235\over208k}}(\log
 x)^{{18627\over8320}(1-1/k)}} $.

\vbox{\vskip1true cm}

{\bf 2.~Auxiliary results. } \bsk

{\bf Lemma 1. } (Vaaler's approximation of fractional parts by
trigonometric polynomials.) { \it For arbitrary $w\in\R$ and any
integer $H>1$, let $\psi(w):= w-[w]-\h$, $$ \psi_H(w):=-\sum_{0<h<H}
{\al_{h,H}}\,\sin(2\pi hw)\,,\qquad \psi_H^*(w):=\sum_{0<h<H}
{\be_{h,H}}\,\cos(2\pi hw)\,,
$$ where, for $h=1,\dots,H-1$,
$$ \al_{h,H}:={1\over\pi h}\,\rho\({h\over H}\)\,,\qquad \be_{h,H}:={1\over H}\(1-{h\over
H}\)\,, $$ and $$ \rho(\xi)= \pi\xi(1-\xi)\cot(\pi\xi)+\xi
\qquad\qquad (0<\xi<1)\,. $$ Then the following inequality holds
true:
$$ \abs{\psi(w)-\psi_H(w)} \le \psi_H^*(w) + {1\over2H}\,. $$ } \bsk

{\bf Proof. } This is one of the main results in Vaaler [27]. A very
well readable exposition can also be found in the monograph by
Graham and Kolesnik [3]. \bsk\msk

{\bf Lemma 2.\quad(A) } { \it Let $F\in C^4[A, B]$, and suppose
that, for positive parameters $X, Y, Z$, we have $1\ll B-A\ll X$ and
$$ F^{(j)}\ll X^{2-j} Y^{-1} \for j=2,3,4, \quad
\abs{F''}\ge c_0 Y^{-1}\,,
$$ throughout the interval $[A,B]$, with some constant $c_0>0$.
Let $\J'$ denote the image of $[A,B]$ under $F'$, and $F^*$ the
inverse function of $F'$. Set further $$ r(\xi):=\cases{0 & if
$F'(\xi)\in\Z$,\cr \min\(\norm{F'(\xi)}^{-1},\sqrt{Y}\) &
else,\cr}$$ where $\norm{\cdot}$ stands for the distance from the nearest
integer. Then, with $e(w)=e^{2\pi iw}$ as usual,
$$ \eq{\sum_{A<n\le B} e(F(n)) =&\ e\({{\rm sgn}(F'')\over8}\)
\summe{*}_{\l\in\J'}
{1\over\sqrt{\abs{F''(F^*(\l))}}}\,e(F(F^*(\l))-\l F^*(\l)) + \cr
& + O(r(A)) + O(r(B)) + \O{\log(2+\,{\rm length}(\J'))}\,, \cr }$$
with the notation \par \vbox{$$\summe{*}_{a\le n\le b}\Phi(n) =
\h(\chi_\Zi(a)\Phi(a)+\chi_\Zi(b)\Phi(b)) + \sum_{a<n<b}
\Phi(n)\,,$$ where $\chi_\Zi$ is the indicator function of the
integers.} \msk {\bf(B)\quad} Let further $G\in C^2[A, B]$, with
$G^{(j)}\ll X^{-j} Z$ for $j=0,1,2$, then it follows that
$$ \eq{\sum_{A<n\le B} G(n)\,e(F(n)) =&\ e\({{\rm sgn}(F'')\over8}\) \sum_{\l\in\J'}{G(F^*(\l))
\over\sqrt{\abs{F''(F^*(\l))}}}\,e\(F(F^*(\l))-\l F^*(\l)\)  + \cr
& + \O{Z\(\sqrt{Y}+\log(2+\,{\rm length}(\J'))\)}\,. \cr }
$$ }\bsk

{\bf Proof. } Transformation formulas of this kind are quite common,
though often with worse error terms. For part (A), see Lemma 2.2 in
K\"{u}hleitner and Nowak [22]. Part (B) can be found as f.~(8.47) in the
recent monograph [12] of H.~Iwaniec and E.~Kowalski. \bsk\msk

{\bf The functions $H_{a,b,1}(x)$ and $H_{a,b,2}(x)$. } \msk For
what follows, compare Kr\"{a}tzel [18], and also his monograph [13]. We
start with the definition of {\it generalized Bessel functions} $$
J_\nu^{(\eta)}(x) =
{2\over\sqrt{\pi}\,\Gamma(\nu+1-1/\eta)}\({x\over2}\)^{\eta\nu/2}
\Int_0^1 (1-t^\eta)^{\nu-1/\eta}\cos(xt)\d t\,, $$ for reals
$\eta\ge1$, $\nu>1/\eta$, and real $x>0$. Further, let $$
\psi_\nu^{(\eta)}(x)=
2\sqrt{\pi}\,\Gamma(\nu+1-1/\eta)\sum_{n=1}^\infty \({x\over\pi
n}\)^{\eta\nu/2}J_\nu^{(\eta)}(2\pi nx)\,, $$ the series converging
absolutely for $\nu>1/\eta$. On the basis of these functions, we
define for reals $a\ge b\ge2$, $$ H_{a,b,1}(x) =
{2\Gamma^2(1/b)\over b\Gamma(2/b)} \psi_{3/a}^{(a)}(x)\,,
\eqno(2.1)$$ and
$$ H_{a,b,2}(x) = 8x\Int_0^1
t^{a-1}(1-t^a)^{1/a-1}\,\psi_{2/b}^{(b)}(xt)\d t\,. \eqno(2.2)$$ Of
course it is important to know the asymptotic behavior of these
functions. Since an asymptotics for the generalized Bessel functions
is provided by [13, Lemma 3.11], it readily follows that
$$ \eq{H_{a,b,1}(x) &= C_1(a,b)\,x^{2-2/a}\sum_{n=1}^\infty
{\sin(2\pi nx-\pi/a)\over n^{1+2/a}}\quad +\ O(x)\,, \cr C_1(a,b) &
:= {2\Gamma^2(1/b)\over
b\Gamma(2/b)}\,{2\over\pi}\({a\over2\pi}\)^{2/a}\Gamma\(1+{2\over a
}\) \,.\cr } \eqno(2.3)$$ To expand $H_{a,b,2}(x)$, we approximate
$J_{2/b}^{(b)}(2\pi nxt)$ again by [13, Lemma 3.11], deriving $$
J_{2/b}^{(b)}(2\pi nxt) = {1\over\sqrt{\pi}}\({b\over2\pi
nxt}\)^{1/b} \sin\(2\pi nxt-{\pi\over2b}\) + \O{(nxt)^{-1}}\,, $$ at
least for $xt$ sufficiently large. The part $\Int_0^{1/2}$ in the
definition of $H_{a,b,2}(x)$ is only an $O(1)$. The arising new
integrals
$$ \Int_{1/2}^1 t^{a-1/b} (1-t^a)^{1/a-1} \sin\(2\pi
nxt-{\pi\over2b}\)\d t $$ are dealt with according to E.T.~Copson
[2, p.~24, formula (11.6)]. This finally yields the asymptotics
$$ \eq{H_{a,b,2}(x) &= C_2(a,b)\,x^{2-1/a-1/b}\sum_{n=1}^\infty
n^{-1-1/a-1/b}\,\sin\(2\pi nx-{\pi\over2a}-{\pi\over2b}\)\ +\
O(x)\,, \cr C_2(a,b) & :=
{16\over\pi}\,{a^{1/a}b^{1/b}\over(2\pi)^{1/a+1/b}}\,\Gamma\(1+{1\over
a}\)\Gamma\(1+{1\over b}\) \,.\cr } \eqno(2.4)$$

\bsk\bsk

{\bf 3.~Preparation of the estimate. }  The "cut-into-slices
approach" gives $$ A_{m,k}(x)=\sum_{|n_1|\le x}
L_k\((x^{mk}-|n_1|^{mk})^{1/m}\)\,, \eqno(3.1) $$ where, for real
$W\ge0$,
$$ L_k(W)=\sum_{|n_2|^k+|n_3|^k\le W} 1 $$ is the number of
lattice points in a Lam\'e disc with length parameter $W^{1/k}$.
According to [13, formulae (3.57) and (3.47)], $$ L_k(W) =
a_k\,W^{2/k} +8I_k(W) -8\Delta_k(W) +O(1)\,, \eqno(3.2) $$ where
$a_k={2\Gamma^2(1/k)\over k\Gamma(2/k)}$ is the area of the unit
Lam\'e disc, $$ \eq{I_k(W) &:= \Int_0^{W^{1/k}}
\psi(u)\d\((W-u^k)^{1/k}\)\,, \cr \Delta_k(W) &:= \sum_{(\he
W)^{1/k}<n\le W^{1/k}} \psi\((W-n^k)^{1/k}\)\,.\cr} \eqno(3.3)$$ By
Euler's formula,
$$ \eq{I_k(W) &= \h \sum_{|n|^k\le W}(W-|n|^k)^{1/k} -
\h\Int_{|u|^k\le W} (W-|u|^k)^{1/k}\d u \cr &= \h \sum_{|n|^k\le
W}(W-|n|^k)^{1/k} \ - {a_k\over4}\,W^{2/k}\,.\cr } $$ Combining this
with (3.1) and (3.2), we obtain $$ \eq{A_{m,k}(x) &= -a_k
\sum_{|n_1|\le x} \(x^{mk}-|n_1|^{mk}\)^{2/(mk)}+4\cS(x)\cr & -8
\sum_{|n_1|\le x}\Delta_k\((x^{mk}-|n_1|^{mk})^{1/m}\)\ + O(x)\,,\cr
} $$ with $$ \cS(x) := \sum_{|n_2|^{mk}+|n_3|^{mk}\le x^{mk}}
\((x^{mk}-|n_3|^{mk})^{1/m}-|n_2|^k\)^{1/k}\,. \eqno(3.4)$$
According to [13, Lemma 3.12], $$ \eq{a_k \sum_{|n_1|\le x}
\(x^{mk}-|n_1|^{mk}\)^{2/(mk)} &= {2\Gamma(1+{2\over
mk})\Gamma({1\over mk})\over mk\,\Gamma(1+{3\over mk})}\,a_k\,x^3 +
a_k\,\Psi_{3/(mk)}^{(mk)}(x) \cr & = \vol(\cB_{m,k})\,x^3 +
H_{mk,k,1}(x)\,, \cr } \eqno(3.5)
$$ with an appeal to the definition (2.1). Hence we arrive at $$
\eq{A_{m,k}(x) &= -\vol(\cB_{m,k})\,x^3 - H_{mk,k,1}(x) + 4\cS(x)\cr
& -8 \sum_{|n_1|\le x}\Delta_k\((x^{mk}-|n_1|^{mk})^{1/m}\)\ +
O(x)\,.\cr }\eqno(3.6)
$$ This formula already outlines the strategy for the completion
of the proof: We will have to evaluate $\cS(x)$, and then to
estimate the multiple fractional parts sum.

\bsk\bsk {\bf 4.~Evaluation of the sum $\cS(x)$. } \bsk

{\bf Proposition. } { \it For large real $x$ and fixed real
numbers $m>1$ and $k>2$, the sum $\cS(x)$ defined in $(3.3)$
satisfies the asymptotic formula
$$ \cS(x) =
{1\over2}\vol(\cB_{m,k})x^3+{1\over2}H_{mk,k,1}(x)+{1\over4}H_{mk,k,2}(x)+
 \O{x^{{339\over208}-{235\over208k}}(\log
 x)^{{18627\over8320}(1-1/k)}}\,.$$} \bsk

{\bf Proof. } Obviously, $$ {1\over4}\cS(x) =
\summe{''}_{n_2^{mk}+n_3^{mk}\le x^{mk}\atop
n_2,n_3\ge0}\((x^{mk}-n_3^{mk})^{1/m}-n_2^k\)^{1/k}=
\summeint{''}{t_1}_{(t_1^k+n_2^k)^m+n_3^{mk}\le x^{mk}\atop t_1,
n_2,n_3\ge0}\,,$$ where $\sum^{''}$ means that terms corresponding
to $n_2n_3=0$ get weight $1\over2$, and the term with $n_2=n_3=0$
gets weight $1\over4$. We subdivide this sum into 6 subsums,
according to the relative size of the three variables of
summation, resp., integration, $t_1, n_2,n_3$. Thus
$\cS_{1,2,3}(x)$ comprehends the case $t_1\le n_2\le n_3$, and so
on. Terms with $n_2=n_3$ get weight $1\over2$, which is symbolized
by the notation $\sum^{'}$. Thus, $$ \cS_{1,2,3}=
\summeint{'}{t_1}_{(t_1^k+n_2^k)^m+n_3^{mk}\le x^{mk}\atop 0\le
t_1\le n_2\le n_3}=\cS_{1,2,3}^*(x)+{1\over16}P_{2,3}(x)+O(x)\,,
$$ where $$  \cS_{1,2,3}^*(x) = \summedoint{}{(t_1,t_2)}_{(t_1^k+t_2^k)^m+n_3^{mk}\le x^{mk}\atop 0\le
t_1\le t_2\le n_3} $$ and $P_{2,3}(x)$ has been defined in [18,
formula (9)], with $(k,\kappa)$ instead of $(mk,k)$. It satisfies
the estimate
$$ P_{2,3}(x) \ll x^{{339\over208}-{235\over208k}}(\log
x)^{(18627/8320)(1-1/k)}\,,$$ according to the first formula on top
of page 768 of [18], if one replaces M.~Huxley's bound in [7] by his
more recent slight improvement in [8]. Further, by the Euler
summation formula, $$ \eq{\cS_{1,2,3}^*(x) &=
\drint{(t_1^k+t_2^k)^m+t_3^{mk}\le x^{mk}\atop 0\le t_1\le t_2\le
t_3} \hskip-1.5true cm\d(t_1,t_2,t_3)\quad - \cr &- \dint{0\le
t_1^{mk}\le t_2^{mk}\le x^{mk}-(t_1^k+t_2^k)^m} \hskip-1.5true cm
\psi\((x^{mk}-(t_1^k+t_2^k)^m)^{1/(mk)}\)\d(t_1,t_2)\ + O(x)\cr &=
\drint{(t_1^k+t_2^k)^m+t_3^{mk}\le x^{mk}\atop 0\le t_1\le t_2\le
t_3} \hskip-1.5true cm\d(t_1,t_2,t_3)\quad +\
{1\over16}H_{mk,k,1}(x) +O(x)\,, \cr} $$ by [18, formula (12)] and
the formulae at the bottom of p.~762 of [18]. Thus, altogether, $$
\cS_{1,2,3}(x) = \hskip-0.8true cm
\drint{(t_1^k+t_2^k)^m+t_3^{mk}\le x^{mk}\atop 0\le t_1\le t_2\le
t_3} \hskip-1.5true cm\d(t_1,t_2,t_3)\ +\ {1\over16}H_{mk,k,1}(x) +
\O{x^{{339\over208}-{235\over208k}}(\log
x)^{{18627\over8320}(1-1/k)}}\,. \eqno(4.1) $$ Similarly, $$
\cS_{1,3,2}(x) = \cS_{1,3,2}^*(x) + {1\over16}P_{3,2}(x)+O(x)\,,
$$ where $$  \cS_{1,3,2}^*(x) = \summedoint{}{(t_1,t_3)}_{(t_1^k+n_2^k)^m+t_3^{mk}\le x^{mk}\atop 0\le
t_1\le t_3\le n_2} $$ and $P_{3,2}(x)$ has been defined in [18,
formula (16)], with $(k,\kappa)$ instead of $(mk,k)$. It satisfies
the same estimate as $P_{2,3}(x)$, appealing to [18, Lemma 2].
Consequently, again by Euler's formula,
$$ \cS_{1,3,2}(x) = \hskip-0.8true cm
\drint{(t_1^k+t_2^k)^m+t_3^{mk}\le x^{mk}\atop 0\le t_1\le t_3\le
t_2} \hskip-1.5true cm\d(t_1,t_2,t_3)\ +\ H_{3,2}^*(x) +
\O{x^{{339\over208}-{235\over208k}}(\log
x)^{{18627\over8320}(1-1/k)}}\,, \eqno(4.2) $$ with $$ H_{3,2}^*(x)
= -  \dint{0\le t_1^{k}\le t_3^{k}\le (x^{mk}-t_3^{mk})^{1/m}-t_1^k}
\hskip-1.5true cm
\psi\(((x^{mk}-t_3^{mk})^{1/m}-t_1^k)^{1/k}\)\d(t_1,t_3)\,. $$
Moreover, $$ \cS_{2,1,3}(x) = \cS_{2,1,3}^*(x) + P_{1,3}^*(x) +
O(x)\,, $$ where $$ \cS_{2,1,3}^*(x) =
\summedoint{}{(t_1,t_2)}_{(t_1^k+t_2^k)^m+n_3^{mk}\le x^{mk}\atop
0\le t_2\le t_1\le n_3} $$ and $$ P_{1,3}^*(x) = -
\summeint{}{t_1}_{2^{-m}(x^{mk}-n_3^{mk})\le t_1^{mk}\le
x^{mk}-n_3^{mk}} \psi\(((x^{mk}-n_3^{mk})^{1/m}-t_1^{k})^{1/k}\)\,.
$$ The substitution $$ z=((x^{mk}-n_3^{mk})^{1/m}-t_1^k)^{1/k}\quad
\Leftrightarrow\quad t_1=((x^{mk}-n_3^{mk})^{1/m}-z^k)^{1/k} $$
yields $$ P_{1,3}^*(x) = \summeint{}{z}_{0\le z^{mk}\le
2^{-m}(x^{mk}-n_3^{mk})} \psi(z)
z^{k-1}\((x^{mk}-n_3^{mk})^{1/m}-z^k\)^{-1+1/k}\,.
$$ Hence, integrating by parts, we conclude that $$ P_{1,3}^*(x)
\ll x\,.  $$ Furthermore, essentially repeating an earlier argument,
$$ \eq{\cS_{2,1,3}^*(x) &= \drint{(t_1^k+t_2^k)^m+t_3^{mk}\le
x^{mk}\atop 0\le t_2\le t_1\le t_3} \hskip-1.5true
cm\d(t_1,t_2,t_3)\quad - \cr &- \dint{0\le t_2^{mk}\le t_1^{mk}\le
x^{mk}-(t_1^k+t_2^k)^m} \hskip-1.5true cm
\psi\((x^{mk}-(t_1^k+t_2^k)^m)^{1/(mk)}\)\d(t_1,t_2)\ + O(x)\cr &=
\drint{(t_1^k+t_2^k)^m+t_3^{mk}\le x^{mk}\atop 0\le t_2\le t_1\le
t_3} \hskip-1.5true cm\d(t_1,t_2,t_3)\quad +\
{1\over16}H_{mk,k,1}(x) +O(x)\,. \cr} $$ Hence also $$
\cS_{2,1,3}^*(x) = \drint{(t_1^k+t_2^k)^m+t_3^{mk}\le x^{mk}\atop
0\le t_2\le t_1\le t_3} \hskip-1.5true cm\d(t_1,t_2,t_3)\quad +\
{1\over16}H_{mk,k,1}(x) +O(x)\,. \eqno(4.3) $$ Next, $$
\cS_{2,3,1}(x) = \cS_{2,3,1}^*(x) + P_{3,1}^*(x)+O(x)\,, $$ with
$$ \cS_{2,3,1}^*(x) = \summedoint{}{(t_1,t_2)}_{(t_1^k+t_2^k)^m+n_3^{mk}\le x^{mk}\atop
0\le t_2\le n_3\le t_1} $$ and $$ P_{3,1}^*(x) = -
\summeint{}{t_1}_{(x^{mk}-n_3^{mk})^{1/m}-n_3^k\le t_1^{k}\le
(x^{mk}-n_3^{mk})^{1/m}}
\psi\(((x^{mk}-n_3^{mk})^{1/m}-t_1^{k})^{1/k}\)\,. $$ Again by [18,
Lemma 1], we bound $P_{3,1}^*(x)$ by
$\O{x^{{339\over208}-{235\over208k}}(\log
x)^{{18627\over8320}(1-1/k)}}$ as well. Further, $$ \cS_{2,3,1}^*(x)
= \drint{(t_1^k+t_2^k)^m+t_3^{mk}\le x^{mk}\atop 0\le t_2\le t_3\le
t_1} \hskip-1.5true cm\d(t_1,t_2,t_3)\quad +\ H_{3,1}^*(x) + O(x)\,,
$$ where $$ \eq{& H_{3,1}^*(x) = - \hskip-0.4true cm \dint{0\le
t_2^{mk}\le x^{mk}-(t_1^k+t_2^k)^m\le t_1^{mk}} \hskip-1.4true cm
\psi\((x^{mk}-(t_1^k+t_2^k)^m)^{1/(mk)} \)\d(t_1,t_2) \cr & =
\hskip-1.4true cm \dint{0\le t_2^{k}\le
z^{k}\le(x^{mk}-z^{mk})^{1/m}-t_2^k} \hskip-1.8true cm \psi(z)
z^{mk-1}
(x^{mk}-z^{mk})^{-1+1/m}((x^{mk}-z^{mk})^{1/m}-t_2^k)^{-1+1/k}\d(z,t_2)\,,
\cr } $$ on the basis of the substitution $$ z =
\(x^{mk}-(t_1^k+t_2^k)^m\)^{1/(mk)} \quad\Leftrightarrow\quad t_1
=\((x^{mk}-z^{mk})^{1/m}-t_2^k\)^{1/k}\,. $$ Observing that
$t_2^k\le\h(x^{mk}-z^{mk})^{1/m}$ and $z^{mk}\le\h x^{mk}$, an
integration by parts shows that $H_{3,1}^*(x)\ll x$. Hence,
altogether, $$ \cS_{2,3,1}(x) = \drint{(t_1^k+t_2^k)^m+t_3^{mk}\le
x^{mk}\atop 0\le t_2\le t_3\le t_1} \hskip-1.5true
cm\d(t_1,t_2,t_3)\quad +\ \O{x^{{339\over208}-{235\over208k}}(\log
x)^{{18627\over8320}(1-1/k)}}\,. \eqno(4.4)  $$ Further on, $$
\cS_{3,1,2}(x) = \cS_{3,1,2}^*(x) + P_{1,2}^*(x)+O(x)\,, $$ with
$$ \cS_{3,1,2}^*(x) = \summedoint{}{(t_1,t_3)}_{(t_1^k+n_2^k)^m+t_3^{mk}\le x^{mk}\atop
0\le t_3\le t_1\le n_2} $$ and $$ P_{1,2}^*(x) = -
\summeint{}{t_1}_{x^{mk}-t_1^{mk}\le (t_1^k+n_2^k)^m\le
x^{mk}\atop 0\le t_1\le n_2} \hskip-0.7true cm
\psi\((x^{mk}-(t_1^k+n_2^k)^{m})^{1/(mk)}\)\,.
$$ In order to bound $P_{1,2}^*(x)$, we substitute $$ z=
\(x^{mk}-(t_1^k+n_2^k)^{m}\)^{1/(mk)} \quad\Leftrightarrow\quad
t_1 = \((x^{mk}-z^{mk})^{1/m}-n_2^k\)^{1/k} $$ and obtain $$
P_{1,2}^*(x) = \summeint{}{z}_{(z^k+n_2^k)^{m}\le x^{mk}-z^{mk}\le
2^m n_2^{mk},\ z\ge0} \hskip-1.3true cm \psi(z)
\varphi_{n_2}(z)\,,
$$ where $$ \varphi_{n_2}(z) = z^{mk-1}
(x^{mk}-z^{mk})^{-1+1/m}\((x^{mk}-z^{mk})^{1/m}-n_2^k\)^{-1+1/k}\,.
$$ The condition $(z^k+n_2^k)^{m}\le x^{mk}-z^{mk}$ on the one
hand implies that $z^{mk}\le\h x^{mk}$, hence
$(x^{mk}-z^{mk})^{-1+1/m}\asymp x^{k-mk}$. On the other hand,
$(x^{mk}-z^{mk})^{1/m}\ge z^k+n_2^k$, thus
$((x^{mk}-z^{mk})^{1/m}-n_2^k)^{-1+1/k}\le z^{1-k}$. Therefore,
altogether, $$\varphi_{n_2}(z)\ll x^{k-mk} z^{mk-k}\ll1\,,$$
uniformly in $n_2$. Since $\varphi_{n_2}(z)$ increases
monotonically with $z$, integration by parts readily gives $$
P_{1,2}^*(x) \ll x\,. $$ Further, $$ \cS_{3,1,2}^*(x) =
\drint{(t_1^k+t_2^k)^m+t_3^{mk}\le x^{mk}\atop 0\le t_3\le t_1\le
t_2} \hskip-1.5true cm\d(t_1,t_2,t_3)\quad +\  H_{1,2}^*(x) +
O(x)\,, $$ where $$ H_{1,2}^*(x) = - \hskip-0.5true cm \dint{0\le
t_3^{mk}\le t_1^{mk}\le2^{-m}(x^{mk}-t_3^{mk})} \hskip-1.3true cm
\psi\(((x^{mk}-t_3^{mk})^{1/m}-t_1^k)^{1/k}\) \d(t_1,t_3)\,. $$
Collecting results, we arrive at $$ \cS_{3,1,2}(x) =
\drint{(t_1^k+t_2^k)^m+t_3^{mk}\le x^{mk}\atop 0\le t_3\le t_1\le
t_2} \hskip-1.5true cm\d(t_1,t_2,t_3)\quad +\  H_{1,2}^*(x) +
O(x)\,. \eqno(4.5) $$ Finally, $$ \cS_{3,2,1}(x) =
\cS_{3,2,1}^*(x) + P_{2,1}^*(x)+O(x)\,, $$ with
$$ \cS_{3,2,1}^*(x) = \summedoint{}{(t_1,t_3)}_{(t_1^k+n_2^k)^m+t_3^{mk}\le x^{mk}\atop
0\le t_3\le n_2\le t_1} $$ and $$ P_{2,1}^*(x) = -
\summeint{}{t_1}_{x^{mk}-n_2^{mk}\le (t_1^k+n_2^k)^m\le
x^{mk}\atop 0\le n_2\le t_1} \hskip-0.7true cm
\psi\((x^{mk}-(t_1^k+n_2^k)^{m})^{1/(mk)}\)\,.
$$ Quite the same analysis as used before for $P_{1,2}^*(x)$
applies again and yields $$ P_{2,1}^*(x) \ll x\,. $$ Furthermore,
$$ \cS_{3,2,1}^*(x) = \drint{(t_1^k+t_2^k)^m+t_3^{mk}\le
x^{mk}\atop 0\le t_3\le t_2\le t_1} \hskip-1.5true
cm\d(t_1,t_2,t_3)\quad +\  H_{2,1}^*(x) + O(x)\,, $$ where $$
H_{2,1}^*(x) = - \hskip-0.5true cm \dint{0\le
t_3^{k}\le(x^{mk}-t_3^{mk})^{1/m}-t_1^k\le t_1^{k}} \hskip-1.3true
cm \psi\(((x^{mk}-t_3^{mk})^{1/m}-t_1^k)^{1/k}\) \d(t_1,t_3)\,. $$
We transform this double integral by the substitution $$ z =
\((x^{mk}-t_3^{mk})^{1/m}-t_1^k\)^{1/k} \quad\Leftrightarrow\quad
t_1= \((x^{mk}-t_3^{mk})^{1/m}-z^k\)^{1/k}\,, $$ to obtain $$
H_{2,1}^*(x) =  \dint{0\le t_3^k\le z^k\le
\he(x^{mk}-t_3^{mk})^{1/m}} \hskip-1.4true cm \psi(z) z^{k-1}
\((x^{mk}-t_3^{mk})^{1/m}-z^k\)^{-1+1/k} \d(t_3,z)\,. $$
Integration by parts shows again that $H_{2,1}^*(x)\ll x$, hence
$$ \cS_{3,2,1}(x) = \drint{(t_1^k+t_2^k)^m+t_3^{mk}\le
x^{mk}\atop 0\le t_3\le t_2\le t_1} \hskip-1.5true
cm\d(t_1,t_2,t_3)\quad +\ O(x)\,.\eqno(4.6) $$ Adding up the results
(4.1) - (4.6), we get $$ \eq{{1\over4}\cS(x) &=
{1\over8}\vol(\cB_{m,k})x^3+{1\over8}H_{mk,k,1}(x)+H_{3,2}^*(x)+H_{1,2}^*(x)\cr
& + \O{x^{{339\over208}-{235\over208k}}(\log
x)^{{18627\over8320}(1-1/k)}}\,.\cr }
$$ Finally, $$ H_{3,2}^*(x) + H_{1,2}^*(x) = - \hskip-0.8true cm \dint{\max(2^m t_1^{mk},(t_1^k+t_3^k)^m)
\le x^{mk}-t_3^{mk},\ t_1,t_3\ge0} \hskip-1.7true cm
\psi\(((x^{mk}-t_3^{mk})^{1/m}-t_1^k)^{1/k}\) \d(t_1,t_3)\,. $$
According to [18], formulae (15) and (20), and p.~763-764, this last
double integral equals ${1\over16}H_{mk,k,2}(x)$. This observation
completes the proof of the Proposition. \msk We conclude this
section by the remark that using the Proposition in (3.6) gives
$$ \eq{A_{m,k}(x) &= \vol(\cB_{m,k})x^3 + H_{mk,k,1}(x)+
H_{mk,k,2}(x)\cr &-8 \sum_{|n_1|\le
x}\Delta_k\((x^{mk}-|n_1|^{mk})^{1/m}\) +
\O{x^{{339\over208}-{235\over208k}}(\log
 x)^{{18627\over8320}(1-1/k)}}\,.\cr} \eqno(4.7)$$ \bsk\bsk

{\bf 5.~Approximating the lattice discrepancy of the Lam\'e disc. }
The lattice point discrepancy of the disc $|u_1|^k+|u_2|^k\le W$ is
expressed by (3.2). In fact, for $I_k(W)$ a very precise evaluation
is known, namely $$ I_k(W) =
{1\over\pi}\({k\over2\pi}\)^{1/k}\Gamma\(1+{1\over
k}\)\,W^{1/k-1/k^2}\sum_{j=1}^\infty j^{-1-1/k}\,\sin\(2\pi
jW^{1/k}-{\pi\over2k}\) + O(1)\,. $$ See [13, p.~148]. For
$\Delta_k(W)$ as defined in (3.3), the sharpest upper bound known
reads
$$ \Delta_k(W)=\O{W^{{131\over208\,k}}(\log W)^{18627\over8320}}\,.
\eqno(5.1)$$ This follows by the argument of Kuba [21], if one uses
Huxley's method in its most recent form [8]. Our present target will
be to approximate $\Delta_k(W)$ by trigonometric polynomials, with a
fairly small $O$-term. This result is intended to be useful for
further calculations involving $\Delta_k(W)$, like a summation with
respect to a third dimension, in the sense of (3.1). \ssk To this
end, recall the notation of Lemma 1 and let $H_W$ denote an
arbitrary map from $\{n\in\Z^+:\ \h W<n^k\le W\ \}$ into the
integers exceeding 1, so that $\max(\log H_W)\ll\log W$. Further,
define $\Delta_{k,H_W}(W)$ and $\Delta_{k,H_W}^*(W)$ analogously to
$\Delta_k(W)$, replacing $\psi$ by $\psi_{H_W(n)}$, resp.,
$\psi_{H_W(n)}^*$, for every value of $n$ in (5.1). By Lemma 1,
$$
|\Delta_k(W)-\Delta_{k,H}(W)|\le\Delta_{k,H}^*(W)+ \sum_{(\he
W)^{1/k}<n\le W^{1/k}}{1\over H_W(n)}\,. \eqno(5.2)
$$ Put for
short $\al=\(\al_{h,H_W(n)}\)_{0<h<H_W(n)}$,
$\be=\(\be_{h,H_W(n)}\)_{0<h<H_W(n)}$, and write
$\gamma=\(\gamma_{h,H_W(n)}\)_{0<h<H_W(n)}$ for either $\al$ or
$\be$. In fact, $\Delta_{k,H_W(n)}(W)$ and
$\Delta_{k,H_W(n)}^*(W)$ can be transformed by the same
calculation, considering exponential sums
$$ E_\gamma :=  \sum_{(\he
W)^{1/k}<n\le W^{1/k}}\quad \sum_{0<h<H_W(n)} \gamma_{h,H_W(n)}\,
e\(-h(W-n^k)^{1/k}\)\,. $$ We divide the range of $n$ by a dyadic
sequence $N_0<N_1<\dots<N_J$ defined by the condition
$$\at{{d\over du}\(-h
(W-u^k)^{1/k}\)}{u=N_j}{}=h\,2^j \quad\iff\quad
N_j=W^{1/k}(1+2^{-jq})^{-1/k}\,,\eqno(5.3)$$ with $$ q:={k\over k-1}
$$ for short, throughout what follows. $J$ is chosen as $$
J:=\left[{\log((W/c_0)^{(1-\la)/k})\over q\log2}\right]+1\,,
\eqno(5.4)$$ where $c_0\in[1,b_0]$ is any real, and $0\le\la<1$,
$b_0\ge1$ are arbitrary constants. Hence $J\ll\log W$, and
$W^{1/k}-N_J\asymp W^{\la/k}$. We further impose the condition that
$H_W$ is constant on each subinterval $]N_j,N_{j+1}]$, say equal to
$H_{W,j}$, for $j=0,1,\dots,J-1$. By simple calculations,
$$N_{j+1}-N_j \asymp 2^{-jq}W^{1/k}\,,\eqno(5.5)$$ and
$$ {d^2\over du^2}\(-h
(W-u^k)^{1/k}\)=(k-1)hWu^{k-2}(W-u^k)^{1/k-2} \asymp h
2^{jq(1/k-2)}\,W^{-1/k} $$ for $u\in[N_j,N_{j+1}]$. Thus we may
apply Lemma 2, part (A), to each of the subintervals
$[N_j,N_{j+1}]$, with the choice of parameters $
X=2^{-jq}W^{1/k}$, $Y=h^{-1}\, 2^{-jq(1/k-2)}\,W^{1/k}$. This
yields
$$ \eq{\sum_{0<h<H_{W,j}}\gamma_{h,H_{W,j}}&\sum_{N_j<n\le N_{j+1}} e\(-h
(W-n^k)^{1/k}\) = \cr = { W^{1/(2k)}\over\sqrt{k-1}}
\sum_{0<h<H_{W,j}}\gamma_{h,H_{W,j}}\,h
&\summe{*}_{\l=2^j\,h}^{2^{j+1}\,h}{(h\l)^{\k/2-1}\over(h^\k+\l^\k)^{1-1/(2\k)}}\,
e\(-W^{1/k}(h^\k+\l^\k)^{1/\k}+\text{1\over8}\)\cr & +\ \O{(\log
W)^2}\,,\cr} \eqno(5.6)$$ since, by the definitions in Lemma 1,
$\gamma_{h,H_{W,j}}\ll h^{-1}$. Using the real and imaginary parts
and taking into account the appropriate portion of (5.2), we readily
infer the following result. \bsk

{\bf Theorem 2. } {\it  Let $W$ be a large parameter, $k>2$ and
$0\le\la<1$ fixed real numbers, $H_{W,j}>1$ integers,
$j=0,1,\dots,J-1$, as explained above, with $\max_j\log
H_{W,j}\ll\log W$, and $J$ given by $(5.5)$. For $\al$ and $\be$
as defined, and $\gamma$ denoting either $\al$ or $\be$, put
$$ \eq{{\textstyle\sum_{H,W,k,j}^{(\gamma)}}:= {W^{1/(2k)}\over\sqrt{k-1}}
\sum_{0<h<H_{W,j}}\gamma_{h,H_{W,j}}\,h & \summe{*}_{\l=2^j
h}^{2^{j+1} h}(h\l)^{\k/2-1}(h^\k+\l^\k)^{-1+1/(2\k)}\,\times\cr
\times&\,
e\(W^{1/k}(h^\k+\l^\k)^{1/\k}-{\textstyle{1\over8}}\)\,.\cr}
$$ Then the remainder term $\Delta_k(W)$ defined in $(3.3)$ for the number
of lattice points in a Lam\'e disc  satisfies the estimate $$
\eq{&\abs{\Delta_{k}(W)-\sum_{j=0}^{J-1}\Im\({\textstyle\sum_{H,W,k,j}^{(\al)}}\)}\cr
& \le \sum_{j=0}^{J-1}\Re\({\textstyle\sum_{H,W,k,j}^{(\be)}}\) +
\sum_{j=0}^{J-1} {\#(\,]N_j,N_{j+1}]\cap\Z)\over H_{W,j}} +
C(W^{\la/k}+(\log W)^3)\,,\cr} $$ with an appropriate constant
$C>0$.} \bsk\msk

{\bf Remark. } This result may be called a "truncated Hardy's
identity for Lam\'e discs". It may well be compared with its
classic counterpart for the circle, namely $$ \sum_{0 \le n \le X}
r(n) - \pi X = \frac{1}{\pi}X^{1/4} \sum_{1\le n \le
Y}\frac{r(n)}{n^{3/4}}
           \cos(2 \pi \sqrt{nX}-3\pi/4)
         +\ O(X^{1/2 + \epsilon}\,Y^{-1/2}) + O(Y^{\epsilon})\,, $$
where $X, Y$ are any large reals: see Ivi\'c [10, f.~(1.9)].
Obviously, the formula in \hbox{Theorem 2} is much more complicated.
Apart from the technical smoothing factors $\al$ and $\be$, this is
due to the fact that a Lam\'e curve has points of curvature zero,
and therefore unavoidable.

\bsk\bsk

{\bf 6.~Estimation of the multiple exponential sum. } Going back to
(4.7), it remains to bound $$ \sum_{|n_1|\le
x}\Delta_k\((x^{mk}-|n_1|^{mk})^{1/m}\)\,. \eqno(6.1)$$ Let $U$ be a
further parameter to be chosen later in terms of $x$. For $|n_1|\le
U$, we use (5.1), to get a total error of
$\O{U\,x^{{131\over208}}(\log x)^{18627\over8320}}$. The remaining
range $U<n_1\le x$ will be divided into dyadic subintervals in
various ways: First, we apply \hbox{Theorem 2,} with $\la=0.47$, to
approximate $\Delta_k\((x^{mk}-|n_1|^{mk})^{1/m}\)$. This gives an
overall error term of $O(x^{1.47})$. Next, assuming w.l.o.g. that
$x\over U$ is a power of 2, we define $\cU_r:=]2^{r-1}U,2^rU]$,
$r=1,\dots,R_1$, where $2^{R_1}U={x\over2}$, hence $R_1\ll\log x$.
Further, $\cV_r:=]x(1-2^{-r}),x(1-2^{-r-1})]$, $r=1,\dots,R_2$,
 with $R_2$ such that $2^{-R_2-1}x<10$,
say. Obviously, also $R_2\ll\log x$, and the remaining range gives
only a small error of $O(x)$. After these preparations, and
appealing to Prop.~2, it will obviously suffice to bound all sums
$$ \sum_{n_1\in\cI}
\at{{\textstyle\sum_{H,W,k,j}^{(\gamma)}}}{W=(x^{mk}-n_1^{mk})^{1/m}}{}\,,
$$ where $\cI$ is any one of the intervals $\cU_r$ or $\cV_r$,
and $0\le j\le J-1$. Summation over $r$ and $j$ will be postponed
to the end of the proof. In order to interchange the order of
summation here, we have to ascertain that
$\at{H_{W,j}}{W=(x^{mk}-n_1^{mk})^{1/m}}{}$ and $J$ do not depend
on $n_1\in\cI$. This is easy for $H$: Call the corresponding
values $H_{(\cI,j)}$, and cf.~formulae (6.8) and (6.13) below. It
is a bit more delicate for $J$: Let us observe that there exists a
constant $b_0>1$ such that $$ {\max\{(x^{mk}-u^{mk})^{1/(mk)}:\
u\in\cI\}\over\min\{(x^{mk}-u^{mk})^{1/(mk)}:\ u\in\cI\}} < b_0
$$ for $\cI$ any one of the intervals $\cU_r$ or $\cV_r$.
Recalling (5.5) and the fact that now $W=(x^{mk}-n_1^{mk})^{1/m}$,
we may choose, for every single $n_1\in\cI$,
$$c_0={(x^{mk}-n_1^{mk})^{1/(mk)}\over\min\{(x^{mk}-u^{mk})^{1/(mk)}:\
u\in\cI\}}\,.$$ Then $c_0\le b_0$ throughout, and $J$ is the same
for all $n_1\in\cI$. Therefore,
$$ \eq{&\sum_{n_1\in\cI}
\at{{\textstyle\sum_{H_{(\cI,j)},W,k,j}^{(\gamma)}}}{W=(x^{mk}-n_1^{mk})^{1/m}}{}=\cr
& =
{e(-{\textstyle{1\over8}})\over\sqrt{k-1}}\sum_{0<h<H_{(\cI,j)}}\gamma_{h,H_{(\cI,j)}}\,
h \summe{*}_{\l=2^j h}^{2^{j+1}
h}(h\l)^{\k/2-1}(h^\k+\l^\k)^{-1+1/(2\k)}\,\times\cr
\times&\,\sum_{n_1\in\cI} (x^{mk}-n_1^{mk})^{1/(2mk)}\,
e\((x^{mk}-n_1^{mk})^{1/(mk)}(h^\k+\l^\k)^{1/\k}\)\,.\cr}
\eqno(6.2)
$$ Let $\Sigma(\cI)$ denote the complex conjugate of the innermost
sum here, and put $N=(h^\k+\l^\k)^{1/\k}\asymp\l$ for short. We
transform this sum by Lemma 2, part (B), with
$G(u)=(x^{mk}-u^{mk})^{1/(2mk)}$,
$F(u)=-N(x^{mk}-u^{mk})^{1/(mk)}$. Apart from the range of
summation and the weight factors $G(n)$, this sum is quite similar
to that in (5.6), with $(h,W,k)$ replaced by $(N,x^{mk},mk)$. Thus
Lemma 2, (B), gives $$ \eq{\Sigma(\cI) = {x
N^{q_1-1/2}\over\sqrt{mk-1}} & \sum_{\l_1\in F'(\cI)}
\l_1^{q_1/2-1} \(N^{q_1}+\l_1^{q_1}\)^{1/q_1-3/2}
e\(-x(N^{q_1}+\l_1^{q_1})^{1/q_1}+{\textstyle{1\over8}}\)\cr & +\
\hbox{ error terms }, \cr} \eqno(6.3) $$ with $$ q_1 := {mk\over
mk-1} $$ for short. To evaluate the error terms in (6.3), we have
to distinguish if $\cI$ stands for $\cU_r$ or $\cV_r$. Let us call
the error terms arising finally in (6.2) $\delta(\cU_r,j)$, resp.,
$\delta(\cV_r,j)$. For $u\in\bar{\cU_r}$, it follows that
$G(u)\asymp x^{1/2}$, and $F''(u)\asymp\l x^{1-mk}(2^rU)^{mk-2}$,
hence $$ {G(u)\over\sqrt{|F''(u)|}}\asymp \l^{-1/2}x^{mk/2} (2^r
U)^{1-mk/2} \for u\in\bar{\cU_r}\,. $$ Thus the corresponding
error term in (6.3) reads $$ \O{\l^{-1/2}x^{mk/2} (2^r
U)^{1-mk/2}} + \O{x^{1/2}\log x}\,. $$ Using this in (6.2) and
summing over $\l$ and $h$, we obtain
$$ \delta(\cU_r,j) = \O{x^{mk/2}\,(2^{r}\,U)^{1-mk/2}\,2^{-jq/2}\,\log x} +
\O{H_{(\cU_r,j)}^{1/2}\,2^{-j(q-1)/2}\,x^{1/2}\,\log
x}\,,\eqno(6.4)
$$ recalling from Lemma 1 that $\gamma_{h,H}\,h$ is bounded. Similarly,
for $u\in\bar{\cV_r}$, we have $G(u)\asymp
(2^{-r})^{1/(2mk)}\,x^{1/2}$, and
$F''(u)\asymp\l\,(2^{-r})^{1/(mk)-2}\,x^{-1}$, thus $$
{G(u)\over\sqrt{|F''(u)|}}\asymp 2^{-r}\,\l^{-1/2}\,x \for
u\in\bar{\cV_r}\,. $$ Therefore, the error terms in (6.3) now read
$$  \O{2^{-r}\,\l^{-1/2}\,x} + \O{(2^{-r})^{1/(2mk)}\,x^{1/2}\,\log
x}\,. $$  Using this in (6.2) and summing over $\l$ and $h$, we
arrive at $$\delta(\cV_r,j)= \O{2^{-r}\,2^{-jq/2}\,x\log x} +
\O{(2^{-r})^{1/(2mk)}\,2^{-j(q-1)/2}\,(H_{(\cV_r,j)}\,x)^{1/2}\log
x}\,. \eqno(6.5) $$ It remains to deal with the main term on the
right hand side of (6.3), let us call it $\Sigma^*(\cI)$. \ssk
{\it Case 1. } $\cI=\cU_r$.\quad For all $v\in F'(\bar{\cU_r})$,
we infer that $v\asymp\l(2^rU)^{mk-1}\,x^{1-mk}\ll\l$. Hence,
putting \def\rux{(2^{-s})}
$$ \Sigma_v^*(\cI):= \sum_{\l_1\in F'(\cI),\ \l_1\le v}
e\(-x(N^{q_1}+\l_1^{q_1})^{1/q_1}\)\,, $$ summation by parts readily
gives $$ \Sigma^*(\cU_r)\ll
\l^{1/2-q_1/2}\,x\,(\l\,(2^rU/x)^{mk-1})^{q_1/2-1}\,\max_{v\in\bar{\cU_r}}\abs{\Sigma_v^*(\cU_r)}\,.
\eqno(6.6) $$ We shall estimate $\Sigma_v^*(\cU_r)$ by van der
Corput's "fourth derivative test": Writing
$\phi(w)=-x(N^{q_1}+w^{q_1})^{1/q_1}$, straightforward computation
gives $$ \eq{\phi^{(4)}(w) =&\ -(q_1-1)\,x\,N^{q_1}
w^{q_1-4}(N^{q_1}+w^{q_1})^{1/q_1-4}\times\cr &\times
\(N^{2q_1}(2-q_1)(3-q_1)+(Nw)^{q_1}(7-4q_1)(1+q_1)
+w^{2q_1}(1+q_1)(2+q_1)\)\,.\cr} $$ For $1<q_1\le{7\over4}$, which
is equivalent to $mk\ge{7\over3}$, the last factor has all
nonnegative coefficients. For $w\in F'(\bar{\cU_r})$, it thus
readily follows that
$\phi^{(4)}(w)\asymp\l^{-3}\,x\,(2^rU/x)^{4-3mk}$. Further, as
explained above, length$(F'(\bar{\cU_r}))\ll\l(2^rU/x)^{mk-1}$. In
what follows, it will be convenient to set $2^r\,U/x=2^{-s}$, with
$1\le s=R_1-r\le R_1-1$. According to [13, p.~34, Theorem 2.6],
uniformly in $v\in\bar{\cU_r}$,
$$ \eq{\Sigma_v^*(\cU_r) \ll&\ \l
\rux^{mk-1}\(\l^{-3}x\rux^{4-3mk}\)^{1/14} \cr &+
\(\l\rux^{mk-1}\)^{3/4} \(\l^{-3}x\rux^{4-3mk}\)^{-1/14}\cr \ll &\
\l^{11/14}\,x^{1/14}\,\rux^{(11mk-10)/14} +
\l^{27/28}\,x^{-1/14}\,\rux^{(27mk-29)/28}\,.\cr}
$$ Using this in (6.6), we conclude that
$$ \Sigma^*(\cU_r) \ll \l^{2/7}\,x^{15/14}\,\rux^{2(mk+1)/7} +
\l^{13/28}\,x^{13/14}\,\rux^{(13mk-1)/28}\,. $$ We use this result
in (6.2), and sum over $\l$ and $h$, to arrive at $$
\eq{&\sum_{n_1\in\cU_r}
\at{{\textstyle\sum_{H_{(\cU_r,j)},W,k,j}^{(\gamma)}}}{W=(x^{mk}-n_1^{mk})^{1/m}}{}
\cr &\ll
(2^j)^{(4k-11)/(14k-14)}\,x^{15/14}\,\rux^{2(mk+1)/7}\,H_{(\cU_r,j)}^{11/14}\cr
&+\
(2^j)^{(13k-27)/(28k-28)}\,x^{13/14}\,\rux^{(13mk-1)/28}\,H_{(\cU_r,j)}^{27/28}
+ \delta(\cU_r,j) \,,\cr} \eqno(6.7)$$ where $\delta(\cU_r,j)$ has
been bounded in (6.4). We balance the first term at the right hand
side against the error term
$O(2^{-jq}\,2^{-s}\,x^2\,H_{(\cU_r,j)}^{-1})$ which arises from the
last sum in Theorem 2. (Recall that $\cU_r$ is of length
$\asymp2^rU$ and $]N_j,N_{j+1}]$ is of length $\asymp2^{-jq}x$: cf.
(5.5). This gives $$ H_{(\cU_r,j)} =
\left[(2^{-j})^{(18k-11)/(25k-25)}\,(2^{-s})^{(10-4mk)/25}\,x^{13/25}\right]+1\,.
\eqno(6.8) $$ Using this in (6.7) yields $$ \eq{&\sum_{n_1\in\cU_r}
\at{{\textstyle\sum_{H_{(\cU_r,j)},W,k,j}^{(\gamma)}}}{W=(x^{mk}-n_1^{mk})^{1/m}}{}
\cr &\ll
(2^{-j})^{7k+11\over25(k-1)}\,(2^{-s})^{4mk+15\over25}\,x^{37/25} +
(2^{-j})^{23k+54\over100(k-1)}\,(2^{-s})^{31mk+35\over100}\,x^{143/100}
\cr &+\ 2^{-jq/2}\,(2^r\,U)^{1-mk/2}\,x^{mk/2}\log x +
(2^{-j})^{9k+7\over25(k-1)}\,(2^r\,U/x)^{5-2mk\over25}\,x^{19/25}\log
x\,.\cr} $$ Summing finally over $j$ and $r$, resp., $s$, we arrive
at\fn{*}{Concerning the last term, we write up this conclusion for
$5-2mk<0$. In the contrary case, we obtain a term $x^{19/25+\ep}$
instead which is negligible.}
$$ \eq{&\sum_{U<n_1\le\h x}\Delta_k\((x^{mk}-|n_1|^{mk})^{1/m}\)
\cr &\ll x^{37/25} + U^{1-mk/2}\,x^{mk/2}\,\log x +
U^{-(2mk-5)/25}\,x^{(2mk+14)/25}\log x\,.\cr} \eqno(6.9) $$ It
remains to balance the middle term on the right hand side here
against the error term $\O{U\,x^{{131\over208}}(\log
x)^{18627\over8320}}$ encountered at the beginning of this section.
This leads to $$ U \asymp x^{104mk-131\over104mk}\,(\log x
)^{-{10307\over4160\,mk}}\,. $$ With that choice of $U$, we finally
infer from (6.9) that $$ \eq{&\sum_{|n_1|\le\h
x}\Delta_k\((x^{mk}-|n_1|^{mk})^{1/m}\) \cr &\ll x^{37/25} +
x^{339mk-262\over208mk}(\log x)^{18627mk-20614\over8320mk} +
x^{{2238mk-655\over2600mk}+\ep}\cr &\ll x^{37/25} +
x^{339mk-262\over208mk}(\log x)^{18627mk-20614\over8320mk} \,.\cr }
\eqno(6.10)$$ It remains to verify that $H_{(\cU_r,j)}$, as defined
by (6.8), is $>1$ throughout. In fact, by (5.5), $$ H_{(\cU_r,j)}
\gg (2^{-J})^{(18k-11)/(25k-25)}\,x^{13/25}\gg
x^{{13\over25}-0.53\({18\over25}-{11\over25k}\)}> x^{0.1384}\,, $$
if $mk\ge{5\over2}$. But the case that ${7\over3}\le mk<{5\over2}$
brings in only a factor $$ \gg\
(2^{-s})^{{1\over25}(10-{28\over3})}\gg x^{-2/75}\,. $$ \msk {\it
Case 2. } $\cI=\cV_r$. \quad By a simple calculation, for all $v\in
F'(\bar{\cV_r})$, we have $v\asymp\l\,2^{r/q_1}$. Hence, with a look
back to (6.3), and $ \Sigma_v^*(\cV_r)$ as before, summation by
parts yields $$ \Sigma^*(\cV_r)\ll
2^{-r}\,\l^{-1/2}\,x\,\max_{v\in\bar{\cV_r}}\abs{\Sigma_v^*(\cV_r)}\,.
\eqno(6.11)
$$ We estimate again $\Sigma_v^*(\cV_r)$ by the "fourth derivative
test". Now $$\phi^{(4)}(w)\asymp\l^{q_1}\,w^{-q_1-3}\,x\asymp
2^{-r(4-3/(mk))}\,\l^{-3}\,x$$ for $w\in\bar{\cV_r}$, hence, again
by [13, p.~34, Theorem 6], uniformly in $v\in\bar{\cV_r}$,
$$ \eq{\Sigma_v^*(\cV_r) \ll&\ 2^{r/q_1}\,\l\, \(2^{-r(4-3/(mk))}\l^{-3}\,x\)^{1/14} \cr &+
\(2^{r/q_1}\,\l\)^{3/4} \(2^{-r(4-3/(mk))}\l^{-3}\,x\)^{-1/14}\cr
\ll &\ 2^{r(10mk-11)/(14mk)}\,\l^{11/14}\,x^{1/14}\, +
2^{r(29mk-27)/(28mk)}\,\l^{27/28}\,x^{-1/14}\,.\cr} $$ Combining
this with (6.11), we get $$ \Sigma^*(\cV_r) \ll
2^{-{r(4mk+11)\over14mk}}\,\l^{2/7}\,x^{15/14} +
2^{r(mk-27)\over28mk}\,\l^{13/28}\,x^{13/14}\,.$$ Using this
result in (6.2), and summing over $\l$ and $h$, we obtain $$
\eq{&\sum_{n_1\in\cV_r}
\at{{\textstyle\sum_{H_{(\cV_r,j)},W,k,j}^{(\gamma)}}}{W=(x^{mk}-n_1^{mk})^{1/m}}{}
\cr &\ll
2^{j{4k-11\over14(k-1)}}\,2^{-r{4mk+11\over14mk}}\,H_{(\cV_r,j)}^{11/14}\,x^{15/14}\cr
&+\
2^{j{13k-27\over28(k-1)}}\,2^{r{mk-27\over28mk}}\,H_{(\cV_r,j)}^{27/28}\,x^{13/14}
+ \delta(\cV_r,j) \,,\cr} \eqno(6.12)$$ where $\delta(\cV_r,j)$
has been bounded in (6.5). The error term coming from the last sum
in Theorem 2 now reads
$O(2^{-jq}\,2^{-r}\,x^2\,H_{(\cV_r,j)}^{-1})$, since $\cV_r$ is of
length $\asymp2^{-r}\,x$. We balance this against the first term
on the right hand side of (6.12) - except for the powers of $2^r$,
with respect to which we simply compensate the factor
$2^{r{mk-27\over28mk}}$ in the next-to-last term of (6.12). This
gives
$$ H_{(\cV_r,j)} =
\left[2^{-j{18k-11\over25(k-1)}}\,2^{-r{mk-27\over27mk}}\,x^{13/25}\right]
+ 1 \,. \eqno(6.13) $$ With this choice of $H_{(\cV_r,j)}$, we
infer from (6.12) that $$ \eq{&\sum_{n_1\in\cV_r}
\at{{\textstyle\sum_{H_{(\cV_r,j)},W,k,j}^{(\gamma)}}}{W=(x^{mk}-n_1^{mk})^{1/m}}{}\quad
+ O(2^{-jq}\,2^{-r}\,x^2\,H_{(\cV_r,j)}^{-1})\cr &\ll
(2^{-j})^{7k+11\over25(k-1)}\,\((2^{-r})^{26mk+27\over27mk}+(2^{-r})^{17/54}\)\,x^{37/25}
+ (2^{-j})^{23k+54\over100(k-1)}\,x^{143/100} \cr &+\
2^{-jq/2}\,2^{-r}\, x\log x +
(2^{-j})^{9k+7\over25(k-1)}\,2^{-r/54}\,x^{19/25}\log
x\,.\cr}\eqno(6.14)
$$ Obviously $H_{(\cV_r,j)}>1$ throughout: As before after (6.10), $$
H_{(\cV_r,j)}\gg x^{0.1384}\,2^{-r{mk-27\over27mk}} \gg
x^{0.1384-1/27} \gg x^{1/10}\,. $$ Since on the right hand side of
(6.14) all exponents of 2 are negative, summation over $r$ and $j$
yields $$ \sum_{\h x<|n_1|\le
x}\Delta_k\((x^{mk}-|n_1|^{mk})^{1/m}\) = O(x^{37/25})\,.$$
Together with (6.10) and (4.7), this completes the proof of
Theorem 1.

\vbox{\vskip 1.5true cm}

\klein \parindent=0pt

\cen{\bf References}  \bsk \def\smc{}

[1] Chamizo F, Crist\'obal E, Ubis A (2009) Lattice points in
rational ellipsoids. J.~Math.~Anal.~Appl. {\bf350}: 283-289 \ssk

[2] Copson ET (1965) Asymptotic expansions. Cambridge: Univ.~Press
\ssk

[3] Graham SW, Kolesnik G (1991) Van der Corput's method of
exponential sums. Cambridge: Univ.~Press \ssk

[4] Haberland K (1993) \"{U}ber die Anzahl der Gitterpunkte in
konvexen Gebieten. Preprint FSU Jena (unpublished)\ssk

[5] Heath-Brown DR (1999) Lattice points in the sphere. In:
Gy\"ory et al. (eds.) Number theory in progress, vol.~2, 883-892.
Berlin: de Gruyter \ssk

[6] Hlawka E (1954) \"Uber Integrale auf konvexen K\"orpern I.
Monatsh Math {\bf 54}: 1-36, II, ibid.~{\bf 54}: 81-99 \ssk

[7] Huxley MN (1996) {Area, lattice points, and exponential sums.}
LMS Monographs, New Ser. {\bf 13}, Oxford: University Press   \ssk

[8] Huxley MN (2003) Exponential sums and lattice points III.
Proc.~London Math.~Soc. (3) {\bf87}: 591-609  \ssk

[9] Iosevich A, Sawyer E, Seeger A (2002) Mean square discrepancy
bounds for the number of lattice points in large convex bodies.
J.~Anal.~Math. {\bf87}: 209-230  \ssk

[10] Ivi\'c A (1996) {The Laplace transform of the square in the
circle and divisor problems.} Stud.~Sci. Math. Hung. {\bf32}:
181-205 \ssk

[11] Ivi\'c A, Kr\"atzel E, K\"uhleitner M, Nowak WG (2006)
Lattice points in large regions and related arithmetic functions:
Recent developments in a very classic topic. In: Proceedings
Conf.~on Elementary and Analytic Number Theory ELAZ'04, held in
Mainz, May 24-28, W.~Schwarz and J.~Steuding eds., Franz Steiner
Verlag, pp.~89-128. \ssk

[12] Iwaniec H, Kowalski E (2004) Analytic Number Theory, AMS
Coll.Publ.~53. Providence, R.I. \ssk

[13] Kr\"atzel E (1988) Lattice points. Berlin: VEB Deutscher
Verlag der Wissenschaften \ssk

[14] Kr\"atzel E (1999) Lattice points in super spheres.
Comment.~Math.~Univ.~Carolinae {\bf40}: 373-391  \ssk

[15] Kr\"atzel E (2000) Analytische Funktionen in der
Zahlentheorie. Wiesbaden: Teubner.  \ssk

[16] Kr\"atzel E (2000) Lattice points in three-dimensional large
convex bodies. Math.~Nachr. {\bf212}: 77--90  \ssk

[17] Kr\"atzel E (2002) Lattice points in three-dimensional convex
bodies with points of Gaussian curvature zero at the boundary.
Monatsh.~Math. {\bf137}: 197--211   \ssk

[18] Kr\"atzel E (2002) Lattice points in some special
three-dimensional convex bodies with points of Gaussian curvature
zero at the boundary. Comment.~Math.~Univ.~Carolinae {\bf43}:
755-771  \ssk

[19] Kr\"{a}tzel E, Nowak WG (2007) Eine explizite Absch\"{a}tzung f\"{u}r die
Gitter-Diskrepanz von Rotationsellipsoiden. Monatsh.~Math
{\bf152}: 45-61 \ssk

[20] Kr\"{a}tzel E, Nowak WG (2008) The lattice discrepancy of bodies
bounded by a rotating Lam\'e's curve. Monatsh.~Math {\bf154}:
145-156  \ssk

[21] Kuba G (1993) On sums of two $k$-th powers of numbers in
residue classes II. Abh.~Math.~Sem.~Hamburg {\bf63}: 87-95  \ssk

[22] K\"{u}hleitner M, Nowak WG (2000) The asymptotic behavior of the
mean-square of fractional part sums. Proc.~Edinburgh Math.~Soc.
{\bf43}: 309-323  \ssk

[23] M\"uller W (1999) Lattice points in large convex bodies.
Monatsh.~Math. {\bf128}: 315--330   \ssk

[24] Nowak WG (1986) On the lattice rest of a convex body in
$\Ri^s$, II. Arch.~Math. (Basel) {\bf47}: 232-237 \ssk

[25] Nowak WG (2008) On the lattice discrepancy of bodies of
rotation with boundary points of curvature zero. Arch.~Math.
(Basel) {\bf90}: 181-192  \ssk

[26] Peter M (2002) Lattice points in convex bodies with planar
points on the boundary. Monatsh.~Math. {\bf135}: 37--57 \ssk

[27] Vaaler JD (1985) Some extremal problems in Fourier analysis.
Bull.~Amer.~Math.~Soc. {\bf12}: 183-216  \ssk

\vbox{\vskip 1.2true cm}

\parindent=1.5true cm

 \bsk

\vbox{Ekkehard Kr\"{a}tzel \ssk

Faculty of Mathematics

University of Vienna

Nordbergstra{\ss}e 15

1090 Wien, \"{O}sterreich \bsk

Werner Georg Nowak \ssk

Institute of Mathematics

Department of Integrative Biology

Universit\"at f\"ur Bodenkultur Wien

Gregor Mendel-Stra{\ss}e 33

1180 Wien, \"Osterreich \ssk

E-mail: {\tt \ nowak@boku.ac.at} }

\bye